\documentclass[12pt,reqno]{amsart}
\usepackage{amssymb,delarray}
\usepackage{amsfonts}
\usepackage{epsfig}
\usepackage[all]{xy}

\makeindex{}

\def\edvo{\rule {6pt}{6pt}}

\newtheorem{thm}{Theorem}
\newtheorem{lem}{Lemma}
\newtheorem{cl}{Assertion}

\newtheorem{cor}{Corollary}

\begin{document}
\baselineskip=14.pt plus 2pt 

\title[
]{On the torsion in a group $\bf F/[M,N]$ in the case of combinatorial asphericity of groups $\bf F/M$ and $\bf F/N$}

\author[]{O.V.~Kulikova}

\address{\newline O.V.Kulikova \newline Lomonosov Moscow State University, \newline Moscow Center for Fundamental and Applied Mathematics \newline (Moscow, Russia)}
\email{olga.kulikova@mail.ru}

\dedicatory{}  \thanks{This work was financially supported by the Russian Science Foundation, project no. 22-11-00075.
}

\begin{abstract} {\it Let $F$ be a non-Abelian free group with basis $A$, $M$ and $N$ be the normal closures of sets $R_M$ and $R_N$ of words in the alphabet $A^{\pm 1}$.
As is known, the group $F/[N, N]$ is torsion-free, but, in general,  torsion in $F/[M, N]$ is possible.
  In the paper of Hartley and Kuz'min (1991), it was proved that
  if $R_M=\{v\}$, $R_N=\{w\}$ and words $v$ and $w$ are not a proper power in $F$, then  $F/[M,N]$ is torsion-free.
  In the present paper a sufficient condition for the absence of torsion in $F/[M,N]$ is obtained, which allows to generalize the result of Hartley and Kuz'min to arbitrary words $v$ and $w$. }
\end{abstract}


\maketitle 

 \setcounter{tocdepth}{2}
\def\st{{\sf st}}

\setcounter{section}{-1}

{\bf 1. Introduction.}

Let $F$ be a non-Abelian free group with basis $A$, and let $M$ and $N$ be normal subgroups of $F$, $[M,N]$ be the mutual commutant of the subgroups $M$ and $N$, i.e. the group generated by all commutators $[g, h]=g^{-1}h^{-1}gh\,\, (g\in M, h\in N)$. Denote by $R^F$ the normal closure of a set $R\subseteq F$ in the group $F$.
 Let $M=R_M^F$, $N=R_N^F$, where $R_M, R_N$ are sets of words in the alphabet $A^{\pm 1}$.

It is known that the group $F/[N,N]$ is torsion-free (\cite{higman})
In general, in groups of type $F/[M, N]$, torsion is possible (\cite{ch_gupta, kuzmin1, kuzmin2}).
For example, for the free group $F$ of rank $\geq 4$, there are elements of order 2 in the center of the group $F/[F^{\prime\prime}, F]$ (\cite{ch_gupta}).

In the paper \cite{kuzmin2}, it was proved that the group $F/[M,N]$ is torsion-free if  $R_M=\{v\}$, $R_N=\{w\}$ and words $v$ and $w$ are not a proper power in $F$. Under these conditions on $M$ and $N$, the groups $F/M$ and $F/N$ are one-relator groups. Since $v$ and $w$ are not proper powers in $F$, $F/M$ and $F/N$ are torsion-free (\cite{Karrass_Magnus_Solitar}). Any one-relator group is combinatorially aspherical (\cite{lind_pioneer}, see also \cite{lind,chiswell}). 

In the present paper, it is proved that, regardless of the presence of torsion in $F/M$ and $F/N$, the absence of torsion in $F/[M,N]$ follows from the combinatorial asphericity of the presentations $<A|R_M>$ and $<A|R_N>$ of groups $F/M$ and $F/N$,
 where
 $R_M, R_N$ are sets of cyclically reduced nonempty words such that no word from $R_M (R_N)$ is conjugated to another word from $R_M (R_N)$ or its inverse.
  It is not clear whether the combinatorial asphericity of $<A|R_M\cup R_N>$ and the same conditions on $R_M\cup R_N$ are sufficient for the absence of torsion in $F/[M,N]$. Of course, if the generalized Whitehead problem is true, i.e. any subpresentation of a combinatorially aspherical presentation is combinatorially aspherical, then this is true.

\vspace{3mm}
{\bf 1. Designations.}

$u=v$, words $u,v$ are equal in a free group;

$u=v$ in $G$, elements $u,v$ are equal in a group $G$;

$<z>$, the cyclic subgroup of a group $G$, generated by an element $z$;

$R^*$, the set of all cyclic permutations of words from a set $R$ and their inverses;

$F(X)$, the free group with basis $X$.

We will often denote a word in the alphabet $A^{\pm 1}$ with the same symbol both as an element in the free group $F(A)$, and as an element in the group $G$ with a presentation $<A|R>$.

We will say that a word $w$ is {\it a proper power} (in the free group $F$) if $w=w_0^l$ for some $w_0$, where $l>1$. If $w=w_0^l$ for $l>0$ and $w_0$ is not a proper power in $F$, then the word $w_0$ will be called  {\it the root} of $w$ and the number $l$ will be called {\it the period} of $w$.

\vspace{3mm}

{\bf 2. Combinatorial asphericity. Pictures.}

Let $\mathcal{P}=\langle A|R\rangle$ be a presentation of a group. We will suppose that the presentation $<A|R>$ {\it satisfies the condition} $\bf RC$ (the relator condition), if any element $r$ from $R$ is a cyclically reduced word in the alphabet $A^{\pm 1}$, not equal to the identity in the free group $F(A)$, and no element $r$ from $R$  is a conjugate of any other element from $R$ or its inverse.

 Asphericity in different variants is found in many papers (\cite{chiswell,lind,sieradski,olsh,pride,gener_b_pr,ratcl}, etc.). Terminology in them is not always the same.  Note the paper \cite{chiswell}, in which different asphericities are compared. 

  Taking into account the condition $\bf RC$, here combinatorial asphericity ($\bf CA$) is understood in a sense that coincides with almost asphericity in \cite{ratcl}, combinatorial asphericity in \cite{chiswell,pride}, asphericity in \cite{huebschmann} (in these papers one can find the exact definition and basic properties of combinatorial asphericity). Briefly, $\langle A|R\rangle$ is combinatorially aspherical ($\bf CA$), if its (2-dimensional) Cayley complex $C(A;R)$ is aspherical.

   By the criterion given in Corollary 2 of Theorem 2.6 of \cite{pride}, combinatorial asphericity can be defined in terms of pictures (i.e., dual objects to van Kampen diagrams). Since we will use pictures in this paper,  below we will define combinatorial asphericity exactly in these terms, first recalling the definition of pictures (according to \cite{gener_b_pr}, \cite{pride}).

A {\it picture} $\mathbf{P}$ over  $\mathcal{P}$ is a finite collections of pairwise disjoint closed disks $\Delta_1,\ldots, \Delta_m$ in a closed disk $D^2$,
together with a finite number of disjoint {\it arcs} $\alpha_1,\ldots,\alpha_l$ properly embedded in $D^2 -\cup_{i=1}^mint\Delta_i$
(where "$int$"\, denotes interior). Loosely speaking, below the disks $\Delta_1,\ldots,\Delta_m$ will be called {\it vertices} of $\mathbf{P}$. An arc
can be either a simple closed curve having trivial intersection  with
 $\partial D^2\cup\partial\Delta_1\cup\ldots\cup\partial\Delta_m$, or a simple non-closed curve which joins two different points of $\partial D^2\cup\partial\Delta_1\cup\ldots\cup\partial\Delta_m$.
The {\it boundary}  of $\mathbf{P}$ is the circle $\partial D^2$ which will be denoted by $\partial \mathbf{P}$.
The {\it corners} of a vertex $\Delta$ of $\mathbf{P}$ are the closures of the connected components of $\partial\Delta -\cup_{j}\alpha_j$.
The {\it regions} of $\mathbf{P}$ are the closures of the connected components of $D^2-(\cup_i\Delta_i\bigcup\cup_j\alpha_j)$.
The picture $\mathbf{P}$ is {\it spherical} if no arc of $\mathbf{P}$ meets $\partial \mathbf{P}$. In this case, it is often convenient to consider $\partial \mathbf{P}$ as a point, and to consider a sphere instead of $D^2$.

Fix an
orientation of the ambient disk $D^2$, thereby determining a sense of positive rotation (i.g. clockwise). Assume that the vertices and arcs of $\mathbf{P}$ are labeled by elements of $\mathcal{P}$ as follows.

(i) Each arc of $\mathbf{P}$ is equipped with a normal orientation (indicated  by an arrow transverse to the arc), and is labeled by an element of $A$.

(ii) Each vertex $\Delta$ of $\mathbf{P}$ is equipped with a sign $\epsilon(\Delta) = \pm 1$ and is labeled by a relator $r(\Delta)\in R$.

For a corner $c$ of a vertex $\Delta$ of $\mathbf{P}$, $W(c)$ denotes the word in the alphabet $A^{\pm 1}$ obtained by reading in order the (signed) labels on the arcs that are encountered in a walk around $\partial\Delta$ in the positive direction, beginning and ending at an interior point of $c$ (with the understanding that if we cross an arc, labeled $y$ say, in the direction of its normal orientation then we read $y$, whereas  if we cross the arc against the orientation we read $y^{-1}$).
The oriented and labeled picture $\mathbf{P}$ is a {\it picture over} $\mathcal{P}$ if for each corner $c$ of each vertex $\Delta$ of $\mathbf{P}$, $W(c)$ is identically equal to a  cyclic permutation of $r(\Delta)^{\epsilon(\Delta)}$.
We call  $W(c)$ the {\it label} of $\Delta$.

A corner $c$ is a {\it basic corner} of $\Delta$ of $\mathbf{P}$ if $W(c)$ is identically equal to $r(\Delta)^{\epsilon(\Delta)}$.  The vertex $\Delta$ has exactly $p$ basic corners, where $p$ is the period of $r(\Delta)$.

A picture $\mathbf{P}$ over $\mathcal{P}$ becomes a {\it based} picture over $\mathcal{P}$ when it is equipped with basepoints as follows.

\begin{itemize}
\item Each vertex $\Delta$ has one {\it basepoint}, which is a selected point in the interior of a basic corner of $\Delta$.

\item $\mathbf{P}$ has a {\it global basepoint}, which is a selected point in $\partial \mathbf{P}$ that does not lie on any arc of $\mathbf{P}$.

\end{itemize}

The {\it boundary label} on a based picture $\mathbf{P}$ over $\mathcal{P}$ is the word 
obtained by reading in order the (signed) labels on the arcs of $\mathbf{P}$ that are encountered in a walk around $\partial \mathbf{P}$ in the positive direction, beginning and ending at the global basepoint. Alteration of the global basepoint or of the orientation of the ambient disk $D^2$ changes the boundary label of $\mathbf{P}$ only up to cyclic permutation and inversion.

There is the following pictorial version of the "van Kampen lemma"\, (\cite{lind}).

\begin{lem}\label{lemVKampen} A word $w$ in the alphabet $A^{\pm 1}$ represents the identity of the group $G$ defined by $\mathcal{P} = \langle A\mid R\rangle$ if and only if there is a based picture $\mathbf{P}$ over $\mathcal{P}$ with boundary label identically equal to $w$.
\end{lem}

A {\it transverse path $\gamma$} in $\mathbf{P}$ over $\mathcal{P}$ is a path in the closure of $D^2-\bigcup_i \Delta_i$ which intersects the arcs of $\mathbf{P}$ only finitely many times (moreover, if the path intersects an arc then it crossed it, and doesn't just touch it), no endpoints of $\gamma$ touches any arc, and whenever $\gamma$ meets  $\partial\mathcal{P}$ or any $\partial\Delta_i$, it does so only in its endpoints. Since we will only consider transverse paths, we will from now on drop the use of the adjective "transverse".

 The subpicture enclosed by a simple closed path $\gamma$ will be called a {\it spherical subpicture} if $\gamma$ intersects no arc. A spherical (sub)picture will be called {\it empty} if it neither consists of any vertex nor any portion of any arc.

\underline{Operations on pictures.}

Generally below we will not distinguish between pictures which are isotopic.

The following operations ("deformations") can be applied to a based picture $\mathbf{P}$ over $\mathcal{P}$ (\cite{gener_b_pr}).

$BRIDGE:$ (Bridge move) See Figure 1.

\unitlength=1mm
\begin{picture}(120,50)(-8,0)
       \thicklines

       \put(70,25){\vector(1,0){10}}
       \thinlines
       \put(40,40){\line(0,-1){30}}
       \put(50,10){\line(0,1){30}}
       \put(38,25){\vector(1,0){5}}
       \put(52,25){\vector(-1,0){5}}

       \multiput(95,17)(0,23){2}{\line(0,-1){7}}
       \multiput(105,10)(0,23){2}{\line(0,1){7}}
       \put(100,33){\oval(10,10)[b]}
       \put(100,17){\oval(10,10)[t]}
       \put(100,24){\vector(0,-1){5}}
       \put(100,26){\vector(0,1){5}}

       \multiput(35,25)(18,0){2}{\rm $x$}
       \put(96,22){\rm $x$}
       \put(102,26){\rm $x$}
\end{picture}

\begin{center}
Figure 1
\end{center}

$FLOAT:$ Deletion of a closed arc  that separates $D^2$ into two parts, one of which contains the global basepoint of $\mathbf{P}$ and all remaining arcs and vertices of  $\mathbf{P}$ (such a closed arc is called a  {\it floating circle}).

$FLOAT^{-1}:$ (Insertion of a floating circle). The opposite of  $FLOAT$.

A {\it folding pair} is a spherical subpicture that contains exactly two vertices such that
\begin{itemize}
\item these two vertices are labeled by the same relator and have opposite signs;

\item the basepoints of the vertices lie in the same region;

\item each arc in the subpicture has an endpoint on each vertex.

\end{itemize}

$FOLD:$ (Deletion of a folding pair). If there is a simple closed path $\beta$ in $D^2$ such that the part of $\mathbf{P}$ encircled by $\beta$ is a folding pair, then remove that part of $\mathbf{P}$ encircled by $\beta$.

$FOLD^{-1}:$ (Insertion of a folding pair). The opposite of $FOLD$.

Let ${\bf X} = \{\mathbf{P}_{\lambda}\mid \lambda\in \Lambda\}$ be a set of based spherical pictures  over $\mathcal{P}$. By an {\it $\bf X$-picture} we mean either a picture $\mathbf{P}_{\lambda}$ from ${\bf X}$ or its mirror image $-\mathbf{P}_{\lambda}$.

$DELETE({\bf X}):$ (Deletion of an $\bf X$-picture). If there is a simple closed path $\beta$ in $D^2$ such that the part of $\mathbf{P}$ enclosed by $\beta$ is a copy of an  ${\bf X}$-picture, then delete that part of $\mathbf{P}$ enclosed by $\beta$.

$DELETE({\bf X})^{-1}:$ (Insertion of an $\bf X$-picture). The opposite of $DELETE({\bf X})$.

Two based spherical pictures are called {\it ${\bf X}$-equivalent} if one of them can be transformed into the other one (up to planar isotopy) by a finite sequence of operations $BRIDGE$, $FLOAT^{\pm 1}$, $FOLD^{\pm 1}$, $DELETE({\bf X})^{\pm 1}$.

A {\it dipole} in a picture over $\mathcal{P}$ consists of an arc which meets two corners $c_1$, $c_2$ in distinct vertices such that
\begin{itemize}
\item the two vertices are labeled by the same relator and have opposite signs;
\item $c_1$ and $c_2$ lie in the same region of the picture;
\item $W(c_1) = W(c_2)^{-1}$.
\end{itemize}

By a {\it complete dipole} over $\mathcal{P}$, we mean a connected based spherical picture over $\mathcal{P}$ that contains just two vertices, and where each arc of the picture constitutes a dipole. Note that a complete dipole is just a folding pair, in that  the vertex basepoints need not lie in the same region. If the relator that labels the two vertices  of the complete dipole has period one, then a complete dipole is exactly the same as a folding pair.  A complete dipole will be called {\it primitive} if the relator labeling its vertices  has root $y$ and period $l>1$, and there is a path joining the vertex basepoints with label $y^f$, where {\rm gcd} $(f,l)=1$.
It follows from Lemma 2.1 \cite{gener_b_pr} that, modulo primitive dipoles, one need not be concerned with choices of vertex basepoints.

 Thanks to Corollary 2 of Theorem 2.6 \cite{pride} (see Theorem 1.6 (2) \cite{gener_b_pr}), combinatorial asphericity can be defined as follows:
 $\mathcal{P}$ is {\it combinatorially aspherical} ($\bf CA$) if and only if every spherical picture over $\mathcal{P}$ is $\bf X$-equivalent to the empty picture, where $\bf X$ is the collection of primitive dipoles for all relators of $\mathcal{P}$, which are a proper power.
  If $FOLD^{-1}$ and $DELETE({\bf X})^{-1}$ are not used, a presentation $\mathcal{P}$ is {\it diagrammatically aspherical} ($\bf DA$). Diagrammatic asphericity implies combinatorial asphericity.

We will say that a presentation $\langle A|R\rangle$ is {\it diagrammatically reducible} ($\bf DR$), if every spherical picture over  $\langle A|R\rangle$, containing a vertex, contains a dipole.
Diagrammatic reducibility implies diagrammatic asphericity.

Diagrammatic asphericity and diagrammatic reducibility are inherited by subpresentations.
 About combinatorial asphericity this is currently unknown (this question is closely related to the well-known Whitehead problem).

There are different methods for checking asphericities
(\cite{lind, lind2, gersten, sieradski}, etc.). For example, one of the sufficient conditions for combinatorial asphericity is the small cancelation condition $C (6)$.

\vspace{3mm}

{\bf 3. The case of ordinary presentations.}

 In this section let's consider that  $F$ is a non-Abelian free group with basis $A$, $R_1, R_2, R$ are sets of words in the alphebet $A^{\pm 1}$, $N_1=R_1^F, N_2=R_2^F, N=R^F$.

   In \cite{kuzmin2}, in proving that the group $F/[N_1,N_2]$ is torsion-free, if $N_1=\{u_1\}^F$, $N_2=\{u_2\}^F$ and words $u_1$ and $u_2$ are not a proper power in $F$, the following properties were used:

 (I) torsionlessness in $F/N_1\cap N_2$ and in $N_1\cap N_2/[N_1,N_2]$ implies torsionlessness in $F/[N_1,N_2]$ (since $F/[N_1,N_2]$ is an extension of the group $F/N_1\cap N_2$ by the subgroup $N_1\cap N_2/[N_1,N_2]$);

 (II) torsionlessness in $F/N_1$ and $F/N_2$ implies torsionlessness in $F/N_1\cap N_2$ (since $F/N_1\cap N_2$ is a subdirect product of $F/N_1$ and $F/N_2$);

 (III) if $H_2(N_1N_2/N_1)$, $H_2(N_1N_2/N_2)$ are torsion-free, $N_1\cap N_2/[N_1,N_2]$ is torsion-free (by Theorem 1 from \cite{kuzmin2}).

To generalize this result, instead of Properties (I) and (II), we will use the following.

\begin{lem}\label{lem0_kruch3}
If the groups $F/[F,N_1]$, $F/[F,N_2]$, $N_1\cap N_2/[N_1,N_2]$ are torsion-free, the group $F/[N_1,N_2]$ is torsion-free.
 \end{lem}
{\it Proof.}
Since the group $F/[F,N_i]$ is torsion-free, torsion in $F/[N_1,N_2]$ is possible only in the subgroup $[F,N_i]/[N_1,N_2]$ ($i=1,2$). So torsion in $F/[N_1,N_2]$ is possible only in the subgroup $[F,N_1]\cap [F,N_2]/[N_1,N_2]$, which contains in $N_1\cap N_2/[N_1,N_2]$.
\edvo

By Lemma \ref{lem0_kruch3} and Property (III), we obtain the following.

\begin{thm}\label{kruch3}
If the groups $F/[F,N_1]$, $F/[F,N_2]$, $H_2(N_1N_2/N_1)$, $H_2(N_1N_2/N_2)$ are torsion-free, the group $F/[N_1,N_2]$ is torsion-free.
 \end{thm}

\begin{lem}\label{lem_kruch3}
If the presentation $<A|R>$ of the group $G=F/N$ satisfies the conditions $\bf RC$ and $\bf CA$, then

(a) the subgroup $N/[F,N]$ of the group $F/[F,N]$ is a free Abelian group generated by $\{r[F,N]\}_{r\in R}$;

(b) the group $F/[F,N]$ is torsion-free;

(c) the group $H_2(U)$ is a free Abelian for any subgroup $U$ of $G$.
 \end{lem}
{\it Proof.}

(a) The proof is similar to one of Theorem 31.1 \cite{olsh} with the peculiarity that in proving the part 1) of the auxiliary Lemma 31.2 from \cite{olsh} one should use pictures (dual objects to diagrams), the definition of $\bf CA$ (every spherical picture over $<A|R>$ can be transformed into the empty picture by a finite sequence of operations $BRIDGE$, $FLOAT^{\pm 1}$, $FOLD^{\pm 1}$, $DELETE({\bf X})^{\pm 1}$, where $\bf X$ is the collection of primitive dipoles for all relators of $R$, which are a proper power) and that deletion and insertion of complete dipoles do not change the algebraic number of $R$-vertices.

(b) Similar to the proof of Theorem 31.2 \cite{olsh}, consider a word $x$ such that $x\neq 1$ in $F/[F,N]$, $x^l=1$ in $F/[F,N]$, where $l>1$ is the order of $x$ in $F/[F,N]$. Since the subgroup $N/[F,N]$ is free Abelian, $x\notin N/[F,N]$. So we have $x^l=1$ in $F/N$ and $x\neq 1$ in $F/N$. Let $s>1$ be the order of $x$ in $F/N$. Hence $s$ divides $l$, $l=l_1s$. Since $<A|R>$ satisfies the conditions $\bf RC$ and $\bf CA$, by Theorem 3 from \cite{huebschmann}, there is a relator $r=r_0^{q}\in R$ with root $r_0$ and an element $f\in F$ such that $s|q$ and $x=fr_0^{t}f^{-1}$ in $F/N$, $t=q/s$. Therefore, replacing $x$ by its conjugate, we can assume that $x=r_0^tz$ in the group $F/[F,N]$, where $z\in N/[F,N]$. Since $z$ belongs to the center of $F/[F,N]$, we have that in $F/[F,N]$,
$$1=x^l=(r_0^tz)^l=r_0^{tl}z^l=r^{l_1}z^l.$$

Since $r[F,N]$ is a basis element in $N/[F,N]$, $l_1$ is divided by $l$. Hence $s=1$, the contradiction.

(c) By Shapiro's Lemma, $H_2(U)\cong H_2(G, B)$, where $B=\mathbb{Z}G\otimes_{\mathbb{Z}U}\mathbb{Z}$ is a right $\mathbb{Z}G$-module. In the proof of Theorem 1.2 \cite{pride}, one obtains that $H_2(G, B)$ is embedded in $B\otimes_{\mathbb{Z}G}N/[N,N]$. By Proposition 1.2 \cite{chiswell}, the relation module $N/[N,N]$ of $G$ decomposes, as a $\mathbb{Z}G$-module, into a direct sum of cyclic submodules $P_r$, $r\in R$, where $P_r$ is generated by the element $\overline{r}=r[N,N]$ subject to the single relation $\widetilde{r}_0\cdot \overline{r}=\overline{r}$, where $r_0$ is the root of $r$ in the free group $F$, $\widetilde{r}_0$ is its image in the group $G$. So $$B\otimes_{\mathbb{Z}G}N/[N,N]\cong \oplus_{r\in R}(\mathbb{Z}G\otimes_{\mathbb{Z}U}\mathbb{Z})\otimes_{\mathbb{Z}G}(\mathbb{Z}Ge_r/\mathbb{Z}G(1-\widetilde{r}_0)e_r).$$
Each group $(\mathbb{Z}G\otimes_{\mathbb{Z}U}\mathbb{Z})\otimes_{\mathbb{Z}G}(\mathbb{Z}G e_r/\mathbb{Z}G(1-\widetilde{r}_0)e_r)$ is a free Abelian group with basis in one-to-one correspondence with the set of double cosets $Ug<\widetilde{r}_0>$. Therefore $B\otimes_{\mathbb{Z}G}N/[N,N]$ is a free Abelian group, and so its subgroup $H_2(G, B)$ is also free.
\edvo

\begin{cor}\label{cor1_kruch3}
 If the presentations $<A|R_1>$ and $<A|R_2>$ satisfy the conditions $\bf RC$ and $\bf CA$, then the group $F/[N_1,N_2]$ is torsion-free.
 \end{cor}
{\it Proof.}
Since $<A|R_i>$ is $\bf CA$, then $F/[F,N_i]$ is torsion-free by Lemma \ref{lem_kruch3} (b) ($i=1,2$), and $H_2(N_1N_2/N_1)$ and $H_2(N_1N_2/N_2)$
are torsion-free by Lemma \ref{lem_kruch3} (c). Hence, by Theorem  \ref{kruch3}, the group $F/[N_1,N_2]$ is torsion-free.
\edvo

Recall that there is no torsion in the group $G=F/N$ with combinatorial aspherical presentation $<A|R>$ if and only if there is no proper power in $R$ (by Theorem 3 and Proposition 1 from \cite{huebschmann}).

Since any one-relator group is combinatorial aspherical (\cite{lind_pioneer}, see also \cite{lind,chiswell}), by Theorem \ref{kruch3}, we obtain the generalization of Corollary 9 from \cite{kuzmin2} to arbitrary words, not necessarily not being a proper power.

 \begin{cor}\label{cor2_kruch3}
 Let $u_1$ and $u_2$ be arbitrary words in the alphabet $A^{\pm 1}$, $R_1=\{u_1\}, R_2=\{u_2\}$.
 Then the group $F/[N_1,N_2]$ is torsion-free.
 \end{cor}

 Since diagrammatic asphericity implies combinatorial asphericity and diagrammatic asphericity is inherited by subpresentations, by Corollary \ref{cor1_kruch3}, we have the following.

  \begin{cor}\label{cor3_kruch4}
If the presentation $<A|R_1\cup R_2>$ satisfies the conditions $\bf RC$ and $\bf DA$, then the group $F/[N_1,N_2]$ is torsion-free.
 \end{cor}

\begin{lem}\label{lem_ratcliff}
If the presentation $\mathcal{P}=<A|R_1\cup R_2>$ satisfies the conditions $\bf RC$ and $\bf CA$ and $R_1^*\cap R_2^*=\emptyset$, then $N_1\cap N_2 =[ N_1, N_2]$.
 \end{lem}
{\it Proof.} It follows from Theorem 1 \cite{ratcl} proved by
Guti\'{e}rrez and Ratcliffe, since $\pi_2(\mathcal{P}) = {\bf T}$ and $\eta({\bf T}) = 0$ in the notation  of Theorem 1.3 \cite{pride}. (This Lemma can also be proved using pictures as in \cite{ok} by gluing together two pictures $\mathbf{P}_{R_1}$ and $\mathbf{P}_{R_2}$ over $<A|R_1>$ and $<A|R_2>$ with boundary label from $N_1\cap N_2$ and transforming the resulting spherical picture by the condition $\bf CA$ so that the insertion of complete dipoles with vertex labels from $R_i$ is performed inside $\mathbf{P}_{R_i}$).
\edvo

 By Corollary \ref{cor3_kruch4} and Lemma \ref{lem_ratcliff}, we obtain the following.

  \begin{cor}\label{cor3_kruch3}
If the presentation $<A|R_1\cup R_2>$ satisfies the conditions $\bf RC$ and $\bf DA$ and $R_1^*\cap R_2^*=\emptyset$, then the group $F/N_1\cap N_2$ is torsion-free.
 \end{cor}

Recall that it is still unknown whether the subpresentation of a combinatorially aspherical presentation is combinatorially aspherical. The affirmative answer to this question is the generalized Whitehead conjecture (this is Conjecture 1 of \cite{ivanov}, see also Theorems 4-5 of \cite{ivanov}, establishing the connection between this conjecture and the Whitehead asphericity conjecture).
   It is not clear whether in Corollary \ref{cor3_kruch4} (or in Corollary \ref{cor3_kruch3}) we can require combinatorial asphericity of the presentation $<A|R_1\cup R_2>$ instead of diagrammatic asphericity.
      If there exist $R_1$ and $R_2$ such that the presentation $<A|R_1\cup R_2>$ satisfies the conditions $\bf RC$ and $\bf CA$, but in the group $F/[N_1, N_2]$ (or with $R_1^*\cap R_2^*=\emptyset$ in the group $F/N_1\cap N_2$) there is a torsion, then according to Corollary \ref{cor1_kruch3} this would give a counterexample to Conjecture 1 of \cite{ivanov}.

   Let's consider some examples using  groups from \cite{olsh}. As follows from Theorem 32.1 of \cite{olsh} (see also \cite{ashmanovOlsh}) and Proposition 1.2 of \cite{chiswell},  diagrammatic asphericity of graded presentations from \cite{olsh} is $\bf CA$. So the (graded) presentations constructed in Chapters 6, 8 of \cite{olsh} (for the free Burnside group of large odd exponent, the Tarski monster and others) is $\bf CA$. By Theorem 6 of \cite{ivanov}, every subpresentation of these presentations is also $\bf CA$.

    EXAMPLE 1. Let $|A|=2$, $<A|R_1>$ be the (graded) presentation constructed in Chapter 6 of \cite{olsh} for the free Burnside group with basis $A$ in the variety $\mathfrak{B}_n$ defined by the identity $x^n=1$ for a sufficiently large odd number $n$, and $<A|R_2>$ be the (graded) presentation constructed in Theorem 28.3 of \cite{olsh} for  the 2-generated simple group, all of whose  proper subgroups are infinite cyclic. Thus $R_1^*\cap R_2^*=\emptyset$, since by construction $R_1$ consists of proper powers and $R_2$ does not contain any proper power. Let us show that the presentation $<A|R_1\cup R_2>$ is not $\bf CA$.  Indeed,  $F/[N_1, N_2]$ is torsion-free by Corollary \ref{cor1_kruch3}. But  $F/N_1\cap N_2$ has torsion (since   otherwise
   $N_2\subseteq N_1\cap N_2$, and so $F/N_2$ would be surjectively mapped onto $F/N_1$, which contradicts the fact that in the 2-generated free Burnside group not all proper subgroups are cyclic (Corollary 35.6 of \cite{olsh})). Hence $N_1\cap N_2\neq[N_1, N_2]$. It remains to apply Lemma \ref{lem_ratcliff}.

    EXAMPLE 2. Let $|A|=2$, $<A|R_1\cup R_2>$ be the (graded) presentation of the group from Theorem 28.1 \cite{olsh} (Tarski monster) or the group from Theorem 28.2 \cite{olsh}, where $R_1$ are relators of the first type (i.e. being a proper power), $R_2$ are relators of the second type (i.e. not being a proper power).
    Although the description of the construction of $R_1$ is very similar to the construction for the Burnside group presentation in Chapter 6 of \cite{olsh}, $F/N_1$ turns out not to be periodic.   Indeed, since the system $R_1\cup R_2$ is independent
       (Theorem 26.3 of \cite{olsh}), $N_2$ is not contained in $N_1$. Therefore there exists an element $g\in N_2$ such that $g\notin N_1$. Since $F/[N_1, N_2]$ is torsion-free by Corollary \ref{cor1_kruch3}, and $F/[N_1, N_2] = F/N_1\cap N_2$ by Lemma \ref{lem_ratcliff}, the element $g(N_1\cap N_2)$ has infinite order in $F/N_1\cap N_2$. This element maps to $(gN_1,N_2)\in F/N_1\times F/N_2$, and hence the element $gN_1$ has infinite order in $F/N_1$.

\begin{cl}\label{claim1}
 Let the basis $A$ consist of two non-empty parts $A_1$ and $A_2$ such that $A_1^{\pm 1}\cap A_2^{\pm 1}=\emptyset$. Let $R_1, R_2$ be sets of words in the alphabets $A_1^{\pm 1}$ and $A_2^{\pm 1}$, respectively.
Then the group $F/[N_1,N_2]=F/N_1\cap N_2$ is torsion-free.
 \end{cl}
 {\it Proof.}
Denote by $G_i$ the group with the presentation $<A_i|R_i>$, $i=1,2$.
 Then $F/N_1\cong G_1\ast F(A_2)$ and $F/N_2 \cong F(A_1)\ast G_2$. Since $R_1\subset F(A_1)$ and $R_2\subset F(A_2)$, then $F/[N_1,N_2]=F/N_1\cap N_2$ (\cite{lind=my}).

 Consider an arbitrary word $w\in F$ such that $w^s=1$ in the group $F/N_1\cap N_2$ for some nonzero integer $s$. Then $w^s=1$ in $F/N_i$ for each $i=1,2$. It follows from the torsion theorem for free products (\cite{lind}) that $w=g_iw_ig_i^{-1}$ in $F/N_i$, where $w_i$ is a word in the alphabet $A_i^{\pm 1}$, $g_i\in F$. Hence $g_1w_1g_i^{-1}=g_2w_2g_2^{-1}$ in $F/N_1N_2 \cong G_1\ast G_2$ and it follows from the conjugacy theorem in free products (\cite{lind}) that $w_i=1$ in $F/N_i$ for each $i=1,2$, i.e. $w=1$ in $F/N_1\cap N_2$.
 \edvo

 Note that under the conditions of Assertion \ref{claim1}, torsion is possible both in $F/N_1$, $F/N_2$, and in  $F/[F,N_1]$, $F/[F,N_2]$. For example, if $A=\{a_1,a_2\}$, $R_1=\{a_1^2\}$, then $F/N_1\cong \mathbb{Z}_2*\mathbb{Z}$. And if $A_1=\{a_1,a_2,a_3,a_4\}, A_2=\{a_5\}$, then there exists an element $g$ of order 2 in the center of the group $F(A_1)/[F(A_1), F(A_1)^{\prime\prime}]$ (\cite{ch_gupta}). Consider $R_1=F(A_1)^{\prime\prime}$. Then the image of the element $g$ in $F/[F,N_1]$ has order 2.

A generalization of Theorem 7 from \cite{kuzmin2} also follows from Lemma \ref{lem_kruch3}.

  \begin{thm}\label{kruch4}
If the presentations $<A|R_1>$ and $<A|R_2>$ satisfies the conditions $\bf RC$ and $\bf CA$, then $N_1\cap N_2/[N_1,N_2]$ is a free Abelian group.
 \end{thm}
{\it Proof.}
 By Theorem 1 from \cite{kuzmin2}, there exists the exact sequence
 $$0\rightarrow H_2(N_1N_2/N_1)\oplus H_2(N_1N_2/N_2)\rightarrow N_1\cap N_2/[N_1,N_2]\rightarrow \widetilde{P}\rightarrow 0,$$
 where $\widetilde{P}$ is free Abelian. Since by Lemma \ref{lem_kruch3} (c), $H_2(N_1N_2/N_1), H_2(N_1N_2/N_2)$ are free Abelian, then the group $N_1\cap N_2/[N_1,N_2]$ is free Abelian.
\edvo

\vspace{3mm}

 {\bf 4. The cases of relative and generalized presentations. }

 Since $\bf DR$ implies $\bf DA$, then Corollary \ref{cor3_kruch3} holds also in the case of diagrammatic reducibility. Let's obtain analogues of this statement for relative and generalized presentations from \cite{BogleyPride}.

The article \cite{BogleyPride} considers relative presentations
$\langle H,X;{\bf R} \rangle$, where $H$ is a group, $X$ is a set of adjoined generators, ${\bf R}$ is a set of words in the alphabet $H\cup X\cup X^{-1}$. Each element of ${\bf R}$ is written in the form $x_1^{\varepsilon_1}h_1x_2^{\varepsilon_2}h_2\ldots x_n^{\varepsilon_n}h_n$, where
$x_i\in X, \varepsilon_i=\pm 1, h_i\in H$, and is {\it cyclically reduced} in the sense that if $h_i=1$ and $x_i=x_{i+1}$ (subscripts mod $n$), then $\varepsilon_i=\varepsilon_{i+1}$. The words in $\bf R$ represent elements of the free product $H\ast  F(X)$.
 The group $G$ defined by the relative presentation $\langle H,X;{\bf R} \rangle$ is the quotient of the free product
$H\ast  F(X)$ by the normal closure $\bf N$ of $\bf R$.
For any subset $\bf S\subset R$, denote by ${\bf S}^{*_x}$ the set of all cyclic permutations of elements of ${\bf S}\cup {\bf S}^{-1}$, which begin with an symbol of $X\cup X^{-1}$.
  The relative presentation $\langle H,X;{\bf R}\rangle$ is called {\it orientable}, if for each $r\in {\bf R}$, $\{r\}^{*_x}\cap {\bf R}=\{r\}$; and no element of ${\bf R}$ is a cyclic permutation of its inverse.

The article \cite{BogleyPride} introduces the concept of a picture over relative
presentation and asphericity of relative
presentation. 
An ordinary presentation can be regarded as a relative presentation with $H=\bf 1$. Asphericity of relative presentations with $H=\bf 1$ is diagrammatic reducibility of ordinary presentations.

\begin{cl}\label{relkruch2} Let
$\bf R_1, R_2$ be two non-empty sets of cyclically reduced element of $H\ast F(X)\setminus
H$ such that $r_1\neq r_2$ in $H\ast F(X)$ for any $r_1\in {\bf R}_1^{*_x}$ and any $r_2\in
{\bf R}_2^{*_x}$;
$\bf N_1,N_2$ be the normal closures of $\bf R_1, R_2$ in $H\ast F(X)$.
Suppose that the relative presentation $\langle H,X;{\bf R_1}, {\bf R_2}\rangle$ is orientable and aspherical.
Then if $H$ is torsion-free,
 $H\ast F(X)/[{\bf N_1,N_2}]$ is torsion-free.
  \end{cl}
{\it Proof.}
  Since $\langle H,X;{\bf R_1}, {\bf R_2}\rangle$ is an orientable and aspherical relative presentation, we have $H\ast F(X)/[{\bf N_1,N_2}]=H\ast F(X)/{\bf N_1}\cap {\bf N_2}$ (\cite{ok2}), and $\langle H,X;{\bf R_i}\rangle$ is also
orientable and aspherical, $i=1,2$.

 Consider an arbitrary word $w\in H\ast F(X)$ such that $w^s=1$ in the group $H\ast F(X)/{\bf N_1}\cap {\bf N_2}$ for some non-zero integer $s$. For each $i=1,2$, if $w\neq 1$ in $H\ast F(X)/{\bf N_i}$, then $w$ has a finite order $s_i>1$ in $H\ast F(X)/{\bf N_i}$. Taking into account that $H$ is torsion-free, by Theorem 1.4 \cite{BogleyPride} $w=g_iz_i^{t_i}g_i^{-1}$ in $H\ast F(X)/{\bf N_i}$, where  $t_i=q_i/s_i$, $s_i|q_i$, for some $r_i=z_i^{q_i}\in {\bf R_i}$  with root $z_i$, $g_i\in H\ast F(X)$.

  If $w\neq 1$ in $H\ast F(X)/{\bf N_1}$ and $H\ast F(X)/{\bf N_2}$, then 
  $g_1z_1^{t_1}g_1^{-1}=g_2z_2^{t_2}g_2^{-1}$  in $H\ast F(X)/{\bf N_1N_2}$, which is impossible by Theorem 1.3 of \cite{BogleyPride} (taking into account Serre's result quoted in \cite{huebschmann}).

 If $w\neq 1$ only in one of these groups, say in $H\ast F(X)/{\bf N_2}$, then ${z_2}^{t_2}=1$ in $H\ast F(X)/{\bf N_1N_2}$, which contradicts Corollary 4 of \cite{BogleyPride}.

  Otherwise $w=1$ in $H\ast F(X)/{\bf N_1}\cap {\bf N_2}$. \edvo

\vspace{3mm}

The article \cite{BogleyPride} also considers a generalized presentation
 $\langle H_i (i\in I);U \rangle$, where $H_i$ are non-trivial groups, $U$ is a set of cyclically reduced elements of $H=\ast_{i\in I}H_i$ of free product length at least 2. The group defined by the generalized presentation $\langle H_i (i\in I);U \rangle$ is isomorphic to the quotient of $H$ by the normal closure $K$
of $U$ in $H$. For any subset $S\subset U$, denote by $S^{*}$ the set of cyclic permutations of elements $S\cup S^{-1}$.

Suppose that for $\langle H_i (i\in I);U \rangle$, $\{u\}^{*}\cap U=\{u\}$ for each $u\in U$.

Let $X=\{x_i|i\in I\}$ be collection of new symbols. For each $u\in U$, having the normal form $u=u_1\ldots u_m$ ($u_{\lambda}\in H_{i_{\lambda}}, \lambda=1,\ldots,m)$, denote the element
$$u_{aug}=x_{i_1}u_1x_{i_1}^{-1}1x_{i_2}u_2x_{i_2}^{-1}1\ldots x_{i_m}u_mx_{i_m}^{-1}1$$
of $H\ast F(X)$. (Here $1$ is the identity of $H$.) Define the relative presentation
$$\langle H_i, x_i (i\in I);U_{aug} \rangle,$$
where $U_{aug}=\{u_{aug}:u\in U\}$.

In \cite{BogleyPride}, under the assumption that
no element of $U$ is a cyclic permutation of its inverse and the relative presentation $\langle H_i, x_i (i\in I);U_{aug} \rangle$ is aspherical, V. A. Bogley and S. J.Pride proved statements for $H/K$ similar to those that we used above to obtain Assertion \ref{relkruch2}.

\begin{cl}\label{relrelkruch2} Let

$H=*_{i\in I}H_i$, where $H_i$ are non-trivial groups;

$U_1, U_2$ be two non-empty sets of cyclically reduced elements of $H$ of free product length at least 2 such that $r_1\neq r_2$ in $H$ for any $r_1\in U_1^{*}$ and any $r_2\in
U_2^{*}$, moreover, no element from $U_1\cup U_2$ is a cyclic permutation of its inverse;

$K_1,K_2$ be the normal closures of $U_1, U_2$ in $H$;



If the relative presentation $ \langle H,X;U_{1,aug}, U_{2,aug}\rangle$ is aspherical and any group $H_i$ $(i\in I)$ is torsion-free, then
 $H/[K_1,K_2]$ is torsion-free.
  \end{cl}
  {\it Proof.}  Similar to the proof of Assertion \ref{relkruch2} by replacing Theorems 1.3, 1.4 \cite{BogleyPride} with Theorems 3.7, 3.8 \cite{BogleyPride}, and Corollary 4 of \cite{BogleyPride} with Lemma 3.2 of \cite{BogleyPride}.\edvo


\end{document}